\documentclass[12pt]{article}
\usepackage[]{amsfonts}
\usepackage{amssymb}
\usepackage{amsmath}
\usepackage{dsfont}

\usepackage{amsthm}

\makeatletter
\newtheorem*{rep@theorem}{\rep@title}
\newcommand{\newreptheorem}[2]{%
\newenvironment{rep#1}[1]{%
 \def\rep@title{#2 \ref{##1}}%
 \begin{rep@theorem}}%
 {\end{rep@theorem}}}
\makeatother

\usepackage[T1]{fontenc} 

\def\couleur(#1 #2 #3)
	{
\special{ps::[end]}}

\def\underset#1#2{\mathrel{\mathop{\kern0pt #2}\limits_{#1}}}

\def\overset#1#2{\mathrel{\mathop{\kern0pt #2}\limits^{#1}}}

\def\bx#1{\setbox1=\hbox{\kern3pt{#1}\kern3pt}			
 \dimen1=\ht1 \advance\dimen1 by 3pt \dimen2=\dp1 \advance\dimen2 by 3pt
 \setbox1=\hbox{\vrule height\dimen1 depth\dimen2\box1\vrule}%
 \setbox1=\vbox{\hrule\box1\hrule}%
 \advance\dimen1 by .4pt \ht1=\dimen1
 \advance\dimen2 by .4pt \dp1=\dimen2 \box1\relax}

\def\wbb#1{\kern#1em}
\def\vci{\vrule  width.02em height1.47ex depth-.0ex}		
\def\11{{\rm\wbb{.2}\vci\wbb{-.37}1}}

\newcommand{\C}{\mathbb{C}}
\newcommand{\Di}{\mathcal{D}}
\newcommand{\supp}{\mathrm{supp}}
\newcommand{\loc}{\mathrm{loc}}

\newcommand{\D}{\mathbb{D}}
\newcommand{\N}{\mathbb{N}}

\newcommand{\Ol}{\mathcal{O}}
\renewcommand{\P}{\mathbb{D}^n}

\newcommand{\de}{\partial}
\newcommand{\debar}{\overline{\de}}
\renewcommand{\Box}{ \hfill$\blacksquare$}

\newtheorem{Theorem}{Theorem}[section]
\newreptheorem{Theorem}{Theorem}
\newtheorem{Lemma}[Theorem]{Lemma}
\newtheorem{Definition}[Theorem]{Definition}
\newtheorem{Remark}[Theorem]{Remark}
\newtheorem{Proposition}[Theorem]{Proposition}
\newreptheorem{Proposition}{Proposition}
\newtheorem{Corollary}[Theorem]{Corollary}
\newreptheorem{Corollary}{Corollary}

\title{On $L^{r}$ hypoellipticity of solutions with compact support of the Cauchy-Riemann equation}

\author{{\sc Eric Amar\ and\  Samuele Mongodi}}

\date{}

\begin{document}

\maketitle
\ \par
\section{Introduction.\quad \quad }
\setcounter{equation}{0}
In this paper, we investigate the inhomogeneous Cauchy-Riemann equation
$$\debar u=g$$
when $g$ has compact support and belongs to some $L^{r}$ space. The question is if it is possible to find a solution $u$ with the same properties, namely, compactly supported and in $L^{r}$.

The $L^{r}$ solvability of the Cauchy-Riemann equation has been discussed by Kerzman for smoothly bounded strongly pseudoconvex domains (see \cite{Kerzman70} and \cite{Kerzman71}), by Fornaess and Sibony in $\C$ with weights and in Runge domains in $\C^2$ (see \cite{fornsib}). Other works on the subject are \cite{Krantz76}, \cite{Jouenne00}, \cite{Fischer01}, \cite{Khidr08}, \cite{ChangLee00}, \cite{AbdelkaderKhidr04} and \cite{Li10}.

The problem of controlling the support of the solution is also widely discussed. In one complex variable, the existence of a compactly supported solution in $\C$ is related to the vanishing of some integrals, resemblant of the \emph{moment conditions} which appear in CR geometry:
$$\int_{\C}g(z)z^kdm_1(z)\;.$$
If these integrals vanish for every $k\in\N$, then there exists a function $u$ such that $\de u/\de\bar{z}=g$ and $\supp u\Subset\{|z|<R\}$ for some $R$.

It is not hard to generalize this result to domains like punctured discs, as we do in Lemma \ref{lmm_1var}.

In higher dimension, it is well known that the existence of a compactly supported solution depends on the vanishing of the cohomology with compact supports; $H^{p,q}_c(\Omega)$ vanishes, for $\Omega\subseteq \C^n$ Stein, if $q<n$. For smooth forms, the existence of a solution compactly supported in a sublevel of some strictly plurisubharmonic exhausting function has also been studied widely, beginning from the work of Andreotti and Grauert (\cite{AndreottiGrauert62}).

Some attempts at controlling the support of the solution were made by Landucci, in the case of smoothly bounded strictly pseudoconvex domains (see \cite{Landucci79} and \cite{Landucci80}).

\medskip

We tackle the problem for a very special class of domains, which generalize the punctured disc: we consider the Stein open domain obtained by removing a compex hypersurface from a polydisc $\P$. Given $f\in\mathcal{O}(\overline{\P})$ with $Z=\{f=0\}$, we consider the domain $\P\setminus Z$: the particular structure of these open sets allows us to give a constructive proof of our results. We will state our results in terms of $(0,q)-$forms, the extension to the $(p,q)-$forms being obvious.

First of all, for $(0,1)-$forms, we have the following.
\begin{repProposition}{prp_p1forms}Let $\Omega\subseteq\C^n$ be a Stein domain and $\omega$ a $(0,1)-$form with coefficients in $L^{r}_c(\Omega)$ such that $\debar \omega=0$. Then there exists a unique $f\in L^{r}_c(\Omega)$ such that $\debar f=\omega$, with $\|f\|_r\leq C\|\omega\|_r$, where $C$ depends only on $\Omega$.\end{repProposition}
This result leaves the question open for $q>1$.

Let $\omega$ be a generic $(0,q)-$form and let us write
$$\omega=\sum_{|J|=n-q}\omega_Jd\hat{\bar{z}}_J\;.$$
We are going to work with the forms satisfying the following condition
\begin{equation*}
(\ast)\qquad\debar_{j_{n-q}}\!\cdots\debar_{j_k}\omega_J\in L^{r}(\C^n)\qquad k=1,\ldots, n-q\;,\quad \forall\ |J|=n-q\;.\end{equation*}

In Theorems \ref{teo_sol_n_poly}, \ref{teo_sol_n1_poly} and \ref{teo_sol_gen_poly}, we show that, given $\omega$ a $(0,q)-$form compactly supported in $\C^n$, with $\debar\omega=0$, with $L^r$ coefficients and satisfying $(\ast)$, we can find a $(0,q-1)-$form $\beta\in L^r_c(\C^n)$ such that $\debar\beta=\omega$. 

This result in $\C^n$ easily gives the corollary
\begin{repCorollary}{cor_vanishing}
Let $\omega$ be a $(0,q)-$form with compact support in $\P\setminus Z$ and satisfying conditions $(\ast)$, then, for any $k\in\N$, we can find a $(0,q-1)-$form $\beta\in L^r_c(\P)$ such that $\debar(f^k\beta)=\omega$. Equivalently, we can find a $(0,q-1)-$form $\eta=f^k\beta$ such that $\eta\in L^r_c(\P)$, $\eta$ is $0$ on $Z$ up to order $k$ and $\debar \eta=\omega$.\end{repCorollary}

Moreover, in the case of $(0,n)-$forms, our construction allows us to obtain a slightly better result.

\begin{repTheorem}{teo_vanish}Let $f\in {\mathcal{O}}(\overline{{\mathbb{D}^n}})$ be a 
holomorphic function in a neighbourhood of the closed unit polydisc in 
${\mathbb{C}}^{n}$ and set $Z=f^{-1}(0).$ If $\omega $ is a $(0,n)$-form 
in $L^{r}_{c}({\mathbb{D}^n}\setminus Z),$  then for every $k\in\N$ we can find a $(0,n-1)$-form $\eta \in L^{r}({\mathbb{D}^n})$ 
such that $f^{-k}\eta\in L^{r}(\mathbb{D}^n)$ and all the coefficients of $\eta$ but at most one are in $L^{r}_c(\mathbb{D}^n\setminus Z)$; moreover, $\eta$ is such that $\overline{\partial }\eta =\omega .$\ \par
\end{repTheorem}

\medskip

The starting point of this work was an incisive question asked by G.~Tomassini and the second author to the first author.

\ \par
\section{Notations.}
\setcounter{equation}{0}We denote by ${\mathbb{D}}$ the unit disc in ${\mathbb{C}}$ 
and by ${\mathbb{D}^n}$ its n-fold product, 
the unit polydisc in ${\mathbb{C}}^{n}$ . The projection from ${\mathbb{C}}^{n}$ 
onto the $j-$ th coordinate will be denoted by $\pi _{j}$ .\ \par
The standard Lebesgue measure on ${\mathbb{C}}^{n}$ will be $dm_{n}$ and 
we will denote by $g\ast _{k}h$ the partial convolution in the $k-$ th 
variable: \ \par

\begin{displaymath} 
(g\ast _{k}h)(z_{1},\cdots ,z_{n}):=\int_{{\mathbb{C}}}^{}{g(\cdots ,}z_{k-1},\zeta 
,z_{k+1},\cdots )h(\cdots ,z_{k-1},z_{k}-\zeta ,z_{k+1},\cdots )dm_{1}(\zeta 
).
\end{displaymath} \ \par
\quad \quad If $T$ is a distribution in ${\mathbb{C}}^{n},$ we set $\overline{\partial }_{j}T=\frac{\partial T}{\partial \overline{z}_{j}},\ 
j=1,\cdots ,n.$\ \par
\quad \quad Let $J=(j_{1},...,j_{q}),\ j_{k}=1,...,n,$ then we define 
$\hat z_{J}\in {\mathbb{C}}^{n-q}$ with coordinates in $J$ deleted. For 
instance $\hat z_{k}=(\cdots ,z_{k-1},z_{k+1},\cdots )\in {\mathbb{C}}^{n-1}.$\ 
\par

\section{On the Cauchy transform.\quad \quad }
\setcounter{equation}{0}Given $\varphi \in {\mathcal{D}}({\mathbb{C}}^{n})$ 
a smooth functions with compact support, the functions  \ \par

\begin{displaymath} 
\zeta \rightarrow \varphi (\cdots ,z_{k-1},\zeta ,z_{k+1},\cdots ),
\end{displaymath} \ \par
for $k=1,\cdots ,n,$ are still smooth and with compact support, contained 
in $\pi _{k}({\mathrm{supp}}\varphi ).$ The Cauchy transform of $\varphi $ 
in the $k^{th}$ variable is\ \par

\begin{displaymath} 
G_{k}(\varphi )(z)=\varphi \ast _{k}:=\int_{{\mathbb{C}}}^{}{\frac{\varphi 
(\cdot \cdot \cdot ,z_{k-1},\ \zeta _{k},z_{k+1},...)}{\pi (\zeta _{k}-z_{k})}\,dm_{1}(\zeta 
_{k})}
\end{displaymath} \ \par
and we know that~\cite{AmarMatheron04}\ \par
\begin{Lemma}\label{lmm_lp_est}We have, with the above notations,\ \par
\quad \quad \quad \quad $\overline{\partial }_{k}G_{k}(\varphi )(z)=\varphi (z)\ \ \ \ \ \ \forall 
\ z\in {\mathbb{C}}^{n}$\ \par
and\ \par
\quad \quad \quad \quad $\displaystyle \displaystyle \left\Vert{G_{k}(\varphi )}\right\Vert _{L^{r}}\leq \displaystyle 
\left\Vert{\frac{1}{\pi \zeta _{1}}}\right\Vert _{L^{1}({\mathbb{D}})}{\times}\displaystyle 
\left\Vert{\varphi }\right\Vert _{L^{r}}.$\ \par
\end{Lemma}
\quad \quad So the Cauchy transform extends as a bounded linear operator 
on $\varphi \in L^{r}_{c}({\mathbb{D}^n}).$ Moreover $G_{k}(\varphi )$ is 
holomorphic in $z_{k}$ outside of the support of $\varphi $ considered 
as a function of $z_{k},\ \hat z_{k}$ being fixed.\ \par
\quad \quad Throughout this note, $f$ will be a given function holomorphic 
in a neighbourhood of $\overline{{\mathbb{D}^n}}$ and $Z=Z(f)$ will denote 
its zero locus. 

The set of directions for which there is a complex lines with that direction contained in $Z$ is an analytic subset of $\mathbb{CP}^{n-1}$ of dimension $n-2$; therefore we can find $n$ linearly independent complex directions not lying in it. So, after a linear change of coordinates, 
for every $1\leq k\leq n,$ we can find a number $N_{k}$ such that, given $n-1$ complex numbers 
$a_{j},\ j\in \{1,\cdots ,n\}\backslash \{k\},$ with $\mid a_{j}\mid <1,$ 
the number of solutions of\ \par

\begin{displaymath} 
f(\cdots ,a_{k-1},z_{k},a_{k+1},\cdots )=0
\end{displaymath} \ \par
as an equation in $z_{k},$ is less than $N_{k}+1.$\ \par
\quad \quad Because these solutions are those of an analytic function, there is 
always a parametrization of them by measurable functions: it is an easy application of \cite[Theorem 7.34]{Shapiro09}; let us denote 
these solutions by $\{c_{1,k}(a),...,\ c_{N_{k},k}(a)\}$ where the functions 
$c_{j,k}=c_{j,k}(a)$ are measurable from $\C^{n-1}$ to $\C$.\ \par
\quad \quad Let $\varphi \in L^{r}_{c}({\mathbb{D}^n}\backslash Z)$ and 
fix $\hat z_{k}\in {\mathbb{D}}^{n-1}\ ;$ denote by $S_{\varphi }(\hat z_{k})$ 
its support as a function of $z_{k}$ which depends on $\hat z_{k}.$ Then, 
by compactness, there exists numbers $\delta_1,\ldots, \delta_n$ such that $S_{\varphi}(\hat{z}_k)$ has distance at least $\delta_k$ from $c_{1,k}(a),\ldots, c_{N_k,k}(a)$, for every $a\in\C^{n-1}$, so there are numbers $r_{j,k}=r_{j,k}(\hat z_{k})\geq \delta >0$ 
such that the disc $D(c_{j,k},\ r_{j,k})$ in the $z_{k}$ variable is \textsl{not} 
in $S_{\varphi }(\hat z_{k}).$ 

However, these discs could intersect without coinciding; suppose that the discs
$$D(c_{j_1, k},\ r_{j_1, k}),\ldots, D(c_{j_h,k},\ r_{j_h,k})$$
form a connected component of the union of all the discs for the variable $z_k$, then we can suppose that $r_{j_i, k}=\delta_k$ for $i=1,\ldots, h$. If the discs 
$$D(c_{j_1,k},\ \delta_k/3N_k),\ldots, D(c_{j_h,k},\ \delta_k/3N_k)$$
are disjoint, then we are done, otherwise, let us consider a connected component of their union and let us suppose, wlog, that it coincides with the union. Obviously, the diameter of such a connected component is less that $\delta_k$, therefore a disc centered in one of the centers with radius $\delta_k$ will enclose the whole connected componente and, by the definition of $\delta_k$ will still be in the complement of $S_\varphi$.

Therefore, we can set all the centers equals to one of them (it is not relevant which one) and take $\delta_k$ as a radius. The functions $c_{j,k}$ will still be measurable. The discs will be then either disjoint or coinciding and their radii will be bounded from below by $\delta_k/3N_k$; we set $\delta=\min\{\delta_1/3N_1,\ldots, \delta_n/3N_n\}$.
\ \par
\ \par
As we already said, $\varphi \ast _{k}=G_{k}(\varphi )$ 
is holomorphic for $z_{k}\notin {\mathbb{D}}$ and for $z_{k}\in D(c_{j,k},\ r_{j,k}).$\ 
\par
\quad \quad This will be precised in the next section with the help of 
the following definitions.\ \par
Let $\varphi \in L^{r}_{c}({\mathbb{D}^n}),$ we define\ \par
\quad \quad \quad \quad $\displaystyle [\varphi ]_{k}(l)=\frac{1}{\pi }\displaystyle \int_{{\mathbb{C}}}^{}{\varphi 
(\cdot \cdot \cdot ,z_{k-1},\zeta _{k},z_{k+1},\cdot \cdot \cdot )\zeta 
_{k}^{l}dm_{1}(\zeta _{k})}\ ;$\ \par
let $\varphi \in L^{r}_{c}({\mathbb{D}^n}\backslash Z),$ we define\ \par
\quad \quad \quad \quad $[\varphi ,j]_{k}(l)=\frac{1}{\pi }\int_{{\mathbb{C}}}^{}{\frac{\varphi 
(\cdots ,z_{k-1},\zeta _{k},z_{k+1},\cdots )}{(\zeta _{k}-c_{j,k})^{m+1}}dm_{1}(\zeta 
_{k})}.$\ \par
\quad \quad We have the following lemma linking this with $\overline{\partial }$ 
equation.\ \par

\begin{Lemma} \label{lmm_1var}Let $\varphi \in L^{r}_{c}({\mathbb{D}^n}\backslash Z)$, then 
the following are equivalent :\ \par
\quad \quad (i) $[\varphi ]_{k}(l)=[\varphi ,j]_{k}(l)=0$ for every $l\in\N$ and $1\leq j\leq N_k$\ \par
\quad \quad (ii) $G_{k}(\varphi )\in L^{r}_{c}({\mathbb{D}^n}\backslash Z)\ (\Rightarrow \overline{\partial 
}_{k}G_{k}(\varphi )=\varphi ).$\ \par
\end{Lemma}
\textbf{Proof: }
Without loss of generality, we can set $k=1$; we notice that, by Lemma \ref{lmm_lp_est}, $G_1(\varphi)$ is in $L^{r}(\C^n)$, so $(ii)$ is equivalent to the compactness of its support. Moreover, we remark that $G_1(\varphi)$ has compact support in $\P\setminus Z$ if and only if for almost every $a=(a_2,\ldots, a_n)\in\C^{n-1}$ the function $z\mapsto G_1(\varphi)(z,a_2,\ldots, a_n)$ has compact support in
$$(\P\setminus Z)\cap\{z_2=a_2,\ldots, z_n=a_n\}=\D\setminus\{c_{1,1}(a),\ldots, c_{1,N_1}(a)\}\;.$$
On the other hand, $[\varphi]_{1}(l)$ and $[\varphi,j]_{1}(l)$ vanish if and only if the integrals that define them vanish for almost every $z_2,\ldots, z_n$. So, we are reduced to the $1$ variable case: let then $c_1,\ldots, c_N$ be points in $\D\subset\C$ and $\phi\in L^{r}_c(\D\setminus\{c_1,\ldots, c_N\})$; we set $G(z)=G_1(\varphi)(z)$.

\medskip

If $(ii)$ holds, for any $h\in\Ol(\D\setminus\{c_1,\ldots, c_N\})$ we have
$$\int_{\C}\varphi(z)h(z)dm_1(z)=\int_{\C}\frac{\partial G(z)}{\partial\bar{z}}h(z)dm_1(z)=-\int_{\C}G(z)\frac{\partial h(z)}{\partial\bar{z}}dm_1(z)=0$$
where we have used Stokes' theorem, as $G(z)$ has compact support. The last integral vanishes because $h$ is holomorphic.

\medskip

On the other hand, suppose that $(i)$ holds and let $K=\supp \varphi$.  Consider $r<1$ such that $K\Subset \D_r=\{|z|<r\}$ and take $z$ with $|z|>r$; then
\begin{eqnarray*}G(z)&=&-\frac{1}{z\pi}\int_{K}\varphi(\zeta)\frac{1}{1-\frac{\zeta}{z}}dm_1(z)=-\frac{1}{\pi z}\int_{K}\varphi(\zeta)\sum_{l\geq 0}\frac{\zeta^l}{z^l}dm_1(\zeta)\\ &=&
-\frac{1}{\pi z}\sum_{l\geq0}z^{-l}\int_{K}\varphi(\zeta)\zeta^ldm_1(\zeta)=-\frac{1}{\pi}\sum_{l\geq0}z^{-l-1}[\varphi]_1(l)\;.\end{eqnarray*}
So, $G(z)=0$ if $|z|>r$, therefore $\supp G(z)\Subset\D$.

Moreover, fix $j$, $1\leq j\leq N$; there exists $r_j>0$ such that the closure of $D(c_j, r_j)=\{|z-c_j|<r_j\}$ does not meet $\supp \varphi(z)$. So, if $|z-c_j|<r_j$, we have
\begin{eqnarray*}G(z)&=&\frac{1}{\pi}\int_{K}\varphi(\zeta)\frac{1}{(\zeta-c_j)-(z-c_j)}dm_1(\zeta)=\frac{1}{\pi}\int_{K}\varphi(\zeta)\frac{1}{\zeta-c_j}\frac{1}{1-(z-c_j)/(\zeta-c_j)}dm_1(\zeta)\\&=&
\frac{1}{\pi}\int_{K}\varphi(\zeta)\frac{1}{\zeta-c_j}\sum_{l\geq0}\frac{(z-c_j)^l}{(\zeta-c_j)^l}dm_1(\zeta)=\frac{1}{\pi}\sum_{l\geq 0}(z-c_j)^l[\varphi,j]_1(l)\;.\end{eqnarray*}
Therefore, by hypothesis, $G(z)=0$ if $|z-c_j|<r_j$, so $\supp G(z)\Subset\D\setminus\{c_1,\ldots, c_N\}$. \Box
\ \par
\ \par
Moreover, we have the following relations between the Cauchy transform and the quantities defined above.

\begin{Lemma}\label{lmm_basic0}If $g$ and $h$ are $L^{r}$ functions, compactly supported in $\P$, and $g\star_1\frac{1}{z_1}=h\star_1\frac{1}{z_1}$ for $z_1\not\in\D$, then $[g]_1(k)=[h]_1(k)$ for every $k$.
\par
\end{Lemma}
\noindent{\bf Proof: } If $z_1\not\in\D$, we have
$$g\star_1\frac{1}{z_1}=\int_{\D}g(\zeta_1,\hat{z}_1)\frac{1}{z_1-\zeta_1}dm_1(\zeta_1)=\frac{1}{z_1}\int_{\D}g(\zeta_1,\hat{z}_1)\frac{1}{1-\frac{\zeta_1}{z_1}}dm_1(\zeta_1)=$$
$$\frac{1}{z_1}\sum_{k\geq0}z_1^{-k}\int_{\D}g(\zeta_1,\hat{z}_1)\zeta_1^kdm_1(\zeta_1)=\sum_{k\geq0}[g]_1(k)z_1^{-k-1}\;.$$
A similar expansion holds for $h$, so that
$$h\star_1\frac{1}{z_1}=\sum_{k\geq0}[h]_1(k)z_1^{-k-1}\;.$$
Therefore, given that $(g-h)\star_1\frac{1}{z_1}=0$ for $z_1\not\in\D$, we have $[g]_1(k)=[h]_1(k)$ for every $k$. \Box
\ \par
\begin{Lemma}\label{lmm_basic1}If If $g$ and $h$ are $L^{r}$ functions, compactly supported in $\P$, and there exists $j\geq1$ such that $g\star_1\frac{1}{z_1}=h\star_1\frac{1}{z_1}$ for every $z_1\in D(c_{j,1}(\hat{z}_1), r_{j,1}(\hat{z}_1))$, then $[g,j]_1(k)=[h,j]_1(k)$ for every $k$.
\par
\end{Lemma}
We omit the proof as it can be easily obtained from the previous one.
\ \par
\ \par
Finally, we recall a result about the solution with compact support of the equation $\debar f=\omega$ when $\omega$ is a $(0,1)-$form with compact support. 

\begin{Proposition}\label{prp_p1forms}Let $\Omega\subseteq\C^n$, $n\geq 2$, be a Stein domain and $\omega$ a $(0,1)-$form with coefficients in $L^{r}_c(\Omega)$ such that $\debar \omega=0$. Then there exists a unique $f\in L^{r}_c(\Omega)$ such that $\debar f=\omega$, with $\|f\|_r\leq C\|\omega\|_r$, where $C$ depends only on $\Omega$.\end{Proposition}
\noindent{\bf Proof: }For the proof in the case $\Omega=\C^n$,  see for instance \cite[Chapter III]{Laurent11}. For a generic $\Omega\subset\C^n$, we notice that if $f_1$ and $f_2$ are two compactly supported (distributional) solutions, then the difference $f_1-f_2$ is $\debar-$closed, that is, a holomorphic function, but then $f_1=f_2$. Moreover, by \cite{Serre55}, $H^{0,1}_c(\Omega)=0$, so there exists at least one distributional solution to $\debar T=\omega$, compactly supported in $\Omega$, on the other hand, we know that there is $f\in L^r_c(\C^n)$, solving $\debar f=\omega$, given, as described in $\cite{Laurent11}$, by convolution with the Cauchy kernel.

Therefore we have $T=f$ and the desired estimate follows. \Box

\section{The coronas construction\quad \quad }
\setcounter{equation}{0}Let $\varphi $ be a function in $L^{r}_{c}({\mathbb{D}^n}\backslash Z)$ 
and consider the Cauchy transform $G_{1}(\varphi )(z)$ ; for a.e. $\hat z_{1},\ G_{1}(z)$ 
is holomorphic in $z_{1}$ in the complement of $S(\hat z_{1}).$\ \par
Because $\pi _{1}({\mathrm{supp}}\varphi )$ is compact in ${\mathbb{D}},$ 
there exists $D(0,r)$ containing $S(\hat z_{1})$ ; let $\delta =(1-r)/3$ 
and define the corona\ \par

\begin{displaymath} 
C_{0}=\{z_{1}\in {\mathbb{D}}\ :\ r+\delta <\mid z_{1}\mid <r+2\delta 
\}\Subset {\mathbb{D}}
\end{displaymath} \ \par
and let $A_{0}=m_{1}(C_{0}).$\ \par
\quad \quad In the same way, set $\delta _{j}(\hat z_{1})=r_{j,1}(\hat z_{1})/3$ 
and define \ \par

\begin{displaymath} 
C_{j}(\hat z_{1})=\{z_{1}\in {\mathbb{D}}\ :\ \delta _{j}(\hat z_{1})\leq 
\mid z_{1}-c_{j,1}\mid \leq 2\delta _{j}(\hat z_{1})\}\Subset {\mathbb{D}}
\end{displaymath} \ \par
and set $A_{j}(\hat z_{1})=1/m_{1}(C_{j}(\hat z_{1})).$\ \par
\begin{Definition} The outer corona component of $\varphi $ is the function\ 
\par
\quad \quad \quad \quad 
\begin{displaymath} 
K^{(1)}_{0}(\varphi )(z)=A_{0}{\11}_{C_{0}}(z_{1})z_{1}G_{1}(\varphi )(z)
\end{displaymath} \ \par
and the inner coronas components of $\varphi $ are the functions\ \par
\quad \quad \quad \quad 
\begin{displaymath} 
K^{(1)}_{j}(\varphi )(z)=A_{j}(\hat z_{1}){\11}_{C_{j}(\hat z_{1})}(z_{1})(z_{1}-c_{j,1})G_{1}(\varphi 
)(z).
\end{displaymath} \ \par
\end{Definition}
\begin{Remark} The outer and inner coronas components of $\varphi $ are 
well defined for a.e. $\hat z_{1}$ , because $\varphi (\cdot ,\hat z_{1})$ 
is in $L^{r}({\mathbb{C}})$ and has compact support for a.e. $\hat z_{1}.$ 
We define exactly the same way the quantities $\displaystyle K^{(k)}_{j}(\varphi )(z)$ 
with respect to the variables $z_{k}.$\ \par
\end{Remark}
\begin{Lemma} The operators $K^{(1)}_{m},\ m=0,\cdots ,N_{1},$ are linear and well defined 
from $L^{r}_{c}({\mathbb{C}}^{n})$ to $L^{r}_{c}({\mathbb{C}}^{n}).$\ 
\par
\end{Lemma}
\textbf{ Proof: } As noted before, $K^{(1)}_{m}(\varphi )$ is well defined 
a.e. and it is obviously linear; moreover, it has compact support in ${\mathbb{D}}$ by definition. 
We know that, by Lemma~\ref{lmm_lp_est}, $\left\Vert{G_{1}(\varphi )}\right\Vert _{L^{r}({\mathbb{C}}^{n})}\leq 
M\left\Vert{\varphi }\right\Vert _{L^{r}({\mathbb{C}}^{n})}$ hence we 
have\ \par
\quad \quad \quad \quad $\parallel K^{(1)}_{0}(\varphi )\parallel _{r}\leq A_{0}\left\Vert{{\11}_{C_{0}}G_{1}(\varphi 
)}\right\Vert _{L^{r}}\leq A_{0}M\left\Vert{\varphi }\right\Vert _{L^{r}},$\ 
\par
where $M:=\left\Vert{\frac{1}{\pi z_{1}}}\right\Vert _{L^{1}({\mathbb{D}})}.$\ 
\par
\quad \quad For $j\geq 1,\ A_{j}(\hat z_{1})=1/m_{1}(C_{j}(\hat z_{1})),$ 
but $m_{1}(C_{j}(\hat z_{1})\geq \delta >0$ uniformly in $\hat z_{1}\in {\mathbb{D}^n}_{n-1}$ 
hence $A_{j}(\hat z_{1})\geq \delta ^{-1}<\infty ,$ uniformly in $\hat z_{1}\in {\mathbb{D}^n}_{n-1}.$ 
So we get\ \par
\quad \quad \quad \quad $\parallel K^{(1)}_{m}(\varphi )\parallel _{r}\leq \left\Vert{A_{m}(\cdot 
)}\right\Vert _{L^{\infty }({\mathbb{D}^n}_{n-1})}{\times}\left\Vert{{\11}_{C_{0}}G_{1}(\varphi 
)}\right\Vert _{L^{r}}\leq \delta ^{-1}M\left\Vert{\varphi }\right\Vert 
_{L^{r}}.$\ \par
\quad \quad So for fixed $\hat z_{1}\in {\mathbb{D}^n}_{n-1}\ K^{(1)}_{m}(\varphi )$ 
has compact support in $z_{1}$ and, because it operates only in $z_{1}$ 
and $\varphi $ has compact support in ${\mathbb{C}}^{n},$ then $K^{(1)}_{m}(\varphi )$ 
has compact support in ${\mathbb{C}}^{n}.$ \hfill$\blacksquare$\ \par
\ \par
\begin{Remark} The operator $K^{(1)}_0$ is also bounded from $L^{r}_c$ to $L^{r}_c$, therefore continuous. The operators $K^{(1)}_m$ for $m\geq1$ are not.\end{Remark}\ \par
\ \par
\quad \quad The following results link the quantities $[\varphi ]_{1}(k)$ 
and $[\varphi ,c_{j,1}]_{1}(k)$ with the corresponding ones for $K^{(1)}_{0}(\varphi )$ 
and $K^{(1)}_{j}(\varphi ).$ 

\begin{Lemma} \label{lmm_series_exp}We have
$$K^{(1)}_0(\varphi)(z)=A_0\mathds{1}_{C_0}(z_1)\sum_{k\geq0}[\varphi]_1(k)z_1^{-k}\;,$$
$$K^{(1)}_j(\varphi)(z)=A_j\mathds{1}_{C_j(\hat{z}_1)}(z_1)\sum_{k\geq0}[\varphi, j]_1(k)(z_1-c_{j,1}(\hat{z}_1))^{k+1}\;,$$
the convergence of the series being uniform in $\overline{C_0}$ or $\overline{C_j}$. \par
\end{Lemma}
{\bf Proof: } If $|z_1|>r+\delta$ and $(\zeta_1, \hat{z}_1)\in\supp\varphi$, then
$$\frac{|\zeta_1|}{|z_1|}\leq \frac{r}{r+\delta}<1\;.$$
so, in particular, if $z_1\in C_0$, then $|z_1|>|\zeta_1|$. Therefore, if $z_1\in C_0$, we have
$$G_1(z)=\frac{1}{\pi}\int_\C\varphi(\zeta_1,\hat{z}_1)\frac{1}{z_1-\zeta_1}dm_1(\zeta_1)=\frac{1}{\pi z_1}\int_\C\varphi(\zeta_1,\hat{z}_1)\frac{1}{1-\frac{\zeta_1}{z_1}}dm_1(\zeta_1)=$$
$$=\frac{1}{\pi z_1}\int_{\C}\varphi(\zeta_1,\hat{z}_1)\sum_{k\geq 0}\frac{\zeta_1^k}{z_1^k}dm_1(\zeta_1)=\frac{1}{z_1}\sum [\varphi]_1(k)z_1^{-k}=\sum_{k\geq0}[\varphi_1](k)z_1^{-k-1}$$
So $K^{(1)}_0(\varphi)=A_0\mathds{1}_{C_0}(z_1)\sum_{k\geq 0}[\varphi_1](k)z_1^{-k}$ and the convergence is obviously uniform on $\overline{C_0}$.

On the other hand, if $z_1\in C_j(\hat{z}_1)$ and $(\zeta_1,\hat{z}_1)\in\supp\varphi$, then
$$\frac{|z_1-c_{j,1}(\hat{z}_1)|}{|\zeta_1-c_{j,1}(\hat{z}_1)|}\leq \frac{2}{3}<1\;,$$
so we have that, for $z_1\in C_j(\hat{z}_1)$,
$$G_1(z)=\frac{1}{\pi}\int_\C\frac{\varphi}{z_1-\zeta_1}dm_1(\zeta)=\frac{1}{\pi}\int_\C\frac{\varphi}{(z_1-c_{j,1}(\hat{z}_1))+(c_{j,1}(\hat{z}_1)-\zeta_1)}dm_1(\zeta_1)=$$
$$\frac{1}{\pi}\int_\C\frac{\varphi}{\zeta_1-c_{j,1}(\hat{z}_1)}\frac{1}{1-\frac{z_1-c_{j,1}(\hat{z}_1)}{\zeta_1-c_{j,1}(\hat{z}_1)} }dm_1(\zeta_1)=$$
$$\frac{1}{\pi}\int_C\frac{\varphi}{\zeta_1-c_{j,1}(\hat{z}_1)}\sum_{k\geq 0}\frac{(z_1-c_{j,1}(\hat{z}_1))^k}{(\zeta_1-c_{j,1}(\hat{z}_1))^k}dm_1(\zeta_1)=\sum_{k\geq 0}[\varphi, j]_1(k)(z_1-c_{j,1}(\hat{z}_1))^k\;.$$
So $K^{(1)}_j(\varphi)=A_j(\hat{z}_1)\mathds{1}_{C_j(\hat{z}_1)}(z_1)\sum_{k\geq0}[\varphi, j]_1(k)(z_1-c_{j,1}(\hat{z}_1))^{k+1}$ and the convergence is obviously uniform on $\overline{C_j}$. \Box
\ \par
\ \par
We set 
$$K^{(1)}(\varphi)=\sum_{m=0}^{N_1}K^{(1)}_m(\varphi)\;.$$
\ \par
\begin{Proposition} We have $[K^{(1)}(\varphi )]_{1}=[\varphi ]_{1}$ 
and $[K^{(1)}(\varphi ),j]_{1}=[\varphi ,j]_{1}.$\ \par
\end{Proposition}
\textbf{Proof: } We divide the proof in several steps.\ \par

\noindent{\emph{1. $[K^{(1)}_0(\varphi)]_1(k)=[\varphi]_1(k)$ - }} We calculate
$$H(z)=K^{(1)}_0(\varphi)\star_1\frac{1}{z_1}=\left(A_0\mathds{1}_{C_0}(z_1)\sum_{k\geq0}[\varphi]_1(k)z_1^{-k}\right)\star_1\frac{1}{z_1}=$$
$$A_0\sum_{k\geq0}[\varphi]_1(k)(\mathds{1}_{C_0}(z_1)z_1^{-k}\star_1\frac{1}{z_1})=A_0\sum_{k\geq0}[\varphi]_1(k)\int_{C_0}\frac{\zeta_1^{-k}}{z_1-\zeta_1}dm_1(\zeta_1)\;.$$
If $z_1\not\in\D$, we know that
$$\int_{C_0}\frac{\zeta_1^{-k}}{z_1-\zeta_1}dm_1(\zeta_1)=A_0^{-1}z_1^{-k-1}$$
so
$$H(z)=A_0\sum_{k\geq0}[\varphi]_1(k)(A_0^{-1}z_1^{-k-1})=\sum_{k\geq0}[\varphi]_1(k)z_1^{-k-1}\;.$$
Then we have that, if $z_1\not\in\D$,
$$H(z)=G_1(z)$$
so, by Lemma \ref{lmm_basic0}, $[\varphi]_1(k)=[K^{(1)}_0(\varphi)]_1(k)$.

\smallskip

\noindent{\emph{2. $[K^{(1)}_j(\varphi)]_1(k)=0$ for $j>0$ - }} We calculate
$$H_j(z)=K^{(1)}_j(\varphi)\star_1\frac{1}{z_1}=\left(A_j(\hat{z}_1)\mathds{1}_{C_j(\hat{z}_1)}(z_1)\sum_{k\geq0}[\varphi,j]_1(k)(z_1-c_{j,1}(\hat{z}_1))^{k+1}\right)\star\frac{1}{z_1}=$$
$$A_j(\hat{z}_1)\sum_{k\geq0}[\varphi,j]_1(k)(\mathds{1}_{C_j(\hat{z}_1)}(z_1-c_{j,1}(\hat{z}_1))^{k+1}\star_1\frac{1}{z_1})=$$
$$A_j(\hat{z}_1)\sum{k\geq0}[\varphi,j]_1(k)\int_{C_{j}}\frac{(\zeta_1-c_{j,1}(\hat{z}_1))^{k+1}}{z_1-\zeta_1}dm_1(\zeta_1)\;.$$
If $|z_1-c_{j,1}(\hat{z}_1)|>r_{j,1}(\hat{z}_1)$, then
$$\int_{C_{j}}\frac{(\zeta_1-c_{j,1}(\hat{z}_1))^{k+1}}{z_1-\zeta_1}dm_1(\zeta_1)=0$$
for every $k\geq0$. Therefore $H_j(z)=0$, so by Lemma \ref{lmm_basic0} $0=[K^{(1)}_j(\varphi)]_1(k)$.

\smallskip

\noindent{\emph{3. $[K^{(1)}_j(\varphi),j]_1(k)=[\varphi,j]_1(k)$ for $j>0$ - }} By direct computation, using Lemma \ref{lmm_series_exp}, we have
$$[K^{(1)}_j\varphi,j]_1(l)=A_j(\hat{z}_1)\sum_{k\geq0}[\varphi,j]_1(k)\int_{C_j(\hat{z}_1)}(\zeta_1-c_{j,1}(\hat{z}_1)^{k+1}(\zeta_1-c_{j,1}(\hat{z}_1))^{-l-1}dm_1(\zeta_1)=$$
$$\sum_{k\geq0}[\varphi,j]_1(k)\delta_{k,l}=[\varphi,j]_1(l)\;.$$

\smallskip

\noindent{\emph{4. $[K^{(1)}_m\varphi,j]_1(k)=0$ if $m\neq j$ - }} By step 2, $H_m(z)=0$ if $|z_1-c_{m,1}(\hat{z}_1)|>r_{m,1}(\hat{z}_1)$, so in particular if $z_1\in D(c_{j,1}(\hat{z}_1), r_{j,1}(\hat{z}_1))$ with $j\neq m$, we have $H(z)=0$. By Lemma \ref{lmm_basic1}, it follows that
$$[K^{(1)}_m(\varphi),j]_1(k)=0$$
if $m\neq j$ and $m\neq0$.

If $m=0$, we notice that, if $|z_1|<r$,
$$H(z)=A_0\sum_{k\geq0}[\varphi]_1(k)\int_{C_0}\frac{\zeta_1^{-k}}{z_1-\zeta_1}dm_1(\zeta_1)$$
and
$$\int_{C_0}\frac{\zeta_1^{-k}}{z_1-\zeta_1}dm_1(\zeta_1)=0$$
for every $k$, as $|z_1|<r<|\zeta_1|$. So $H(z)=0$ and by Lemma \ref{lmm_basic0} we have that $[K^{(1)}_0(\varphi),j]_1(k)=0$ for every $k$. \Box
\ \par
\ \par
\ \par
\begin{Corollary} \label{cor_structure}Let $\varphi \in L^{r}_{c}({\mathbb{D}^n}\backslash Z),$ 
there are $\varphi _{1},...,\varphi _{n},$ all in $L^{r}_{c}({\mathbb{D}^n}\backslash Z)$ 
and such that\ \par
\quad \quad \quad \quad $\varphi =\varphi _{1}+\cdot \cdot \cdot +\varphi _{n},\ \forall i<n,\ 
\forall j=1,...,N_{i},\ [\varphi _{i}]_{i}=[\varphi _{i},j]_{i}=0.$\ \par
\end{Corollary}
\textbf{Proof: } We set $\varphi _{1}:=\varphi -K^{(1)}\varphi $ and we notice that $[\varphi _{1}]_{1}=0,\ [\varphi _{1},j]_{1}=0$ 
for every $j=1,\cdots ,N_{1}.$\ \par
\quad \quad Now, we can repeat this procedure replacing $z_{1}$ by $z_{2}$ 
and $\varphi $ by $K^{(1)}(\varphi )$ ; we will apply then the operators 
$K^{(2)}_{m},$ defined with respect to the variable $z_{2},$ with the 
relative coronas.\ \par
We set $\varphi _{2}:=K^{(1)}\varphi -K^{(2)}K^{(1)}\varphi $ with the 
property that $[\varphi _{2}]_{2}=0,\ [\varphi _{2},j]_{2}=0$ for every 
$j=1,\cdots ,N_{2}.$\ \par
\quad \quad Iterating the algorithm we set $\varphi _{n-1}:=K^{(n-2)}\cdot \cdot \cdot K^{(1)}\varphi -K^{(n-1)}\cdot 
\cdot \cdot K^{(1)}\varphi $ and\ \par
\quad \quad \quad \quad $\varphi _{n}:=\varphi -\varphi _{1}-\cdot \cdot \cdot -\varphi _{n-1}.$\ 
\par
By an easy recursion we have 
\begin{displaymath} 
\varphi _{n}=K^{(n-1)}\cdots K^{(1)}\varphi 
\end{displaymath}  with, of course $\varphi =\varphi _{1}+\cdots +\varphi _{n}.$\ 
\par
\quad \quad So finally we find a decomposition $\varphi =\varphi _{1}+\cdots +\varphi _{n}$ 
such that, for $i<n,$ we have $[\varphi _{i}]_{i}=0,\ [\varphi _{i},j]_{i}=0$ 
for every $j=1,\cdots ,N_{i}.$ \hfill$\blacksquare$\ \par

\ \par

We have a first result on solvability of the Cauchy-Riemann equation with some control on the support of the solution.

\begin{Theorem} \label{teo_vanish}Let $f\in {\mathcal{O}}(\overline{{\mathbb{D}^n}})$ be a 
holomorphic function in a neighbourhood of the closed unit polydisc in 
${\mathbb{C}}^{n}$ and set $Z=f^{-1}(0).$ If $\omega $ is a $(0,n)$-form 
in $L^{r}_{c}({\mathbb{D}^n}\setminus Z),$  then for every $k\in\N$ we can find a $(0,n-1)$-form $\eta \in L^{r}({\mathbb{D}^n})$ 
such that $f^{-k}\eta\in L^{r}(\mathbb{D}^n)$ and all the coefficients of $\eta$ but at most one are in $L^{r}_c(\mathbb{D}^n\setminus Z)$; moreover, $\eta$ is such that $\overline{\partial }\eta =\omega .$\ \par
\end{Theorem}
\textbf{Proof: }We write
$$\omega=\phi d\bar{z}_1\wedge \cdots\wedge d\bar{z}_n$$
and apply the result of Corollary \ref{cor_structure} to $f^{-k}\phi$. We get
$$f^{-k}\phi=\phi_1+\ldots+\phi_n$$
and $[\phi_i]_i=[\phi_i,h]_i=0$ for $i=1,\ldots, n-1$ and $h=1,\ldots, N_i$. Therefore, by Lemma \ref{lmm_1var}, the functions
$$F_1=\phi_1\star_1\frac{1}{\pi z_1}\;,\quad\hdots\quad,\ \ F_{n-1}=\phi_{n-1}\star_{n-1}\frac{1}{\pi z_{n-1}}$$
are compactly supported in $\mathbb{D}^n\setminus Z$. However, for
$$F_n=\phi_n\star_n\frac{1}{\pi z_n}$$
we only know that $F_n\in L^{r}(\mathbb{D}^n)$. We note that
$$\overline{\partial}(F_1d\hat{\bar{z}}_1+\ldots+ F_nd\hat{\bar{z}}_n)=f^{-k}\phi\;,$$
therefore we define
$$\eta=f^k(F_1d\hat{\bar{z}}_1+\ldots+ F_nd\hat{\bar{z}}_n)$$
and we have
$$\debar \eta=f^k\debar(F_1d\hat{\bar{z}}_1+\ldots+ F_nd\hat{\bar{z}}_n)=\phi\;.$$
It is easy to see that $\eta$ satisfies all the requests of the theorem.  \hfill$\blacksquare$\ \par

\ \par

\section{Obstructions to a solution with compact support}

\quad \quad Let us define the two quantities which tell us when the last 
term in the decomposition from Corollary \ref{cor_structure} verifies also\ \par
\quad \quad \quad \quad $\forall j=1,...,\ N_{n},\ [\varphi _{n}]_{n}=0,\ [\varphi _{n},j]_{n}=0.$\ 
\par
We note that \ \par

\begin{displaymath} 
\varphi _{n}=K^{(n-1)}\cdots K^{(1)}\varphi 
\end{displaymath} \ \par
and, more precisely, we have\ \par
\quad \quad \quad \quad $\displaystyle \varphi _{n}=\displaystyle \sum_{m_{n-1}=0}^{N_{n-1}}{\cdot \cdot \cdot 
\displaystyle \sum_{m_{1}=0}^{N_{1}}{K^{(n-1)}_{m_{n-1}}\cdot \cdot \cdot 
K^{(1)}_{m_{1}}(\varphi )}}.$\ \par
We set\ \par
\quad \quad \quad \quad $M_{n-1}:=\{(m_{1},...,\ m_{n-1})::m_{j}\leq N_{j}\}\subset {\mathbb{N}}^{n-1}\ 
;$\ \par
\quad \quad \quad \quad $\mu =(m_{1},...,m_{n-1})\in M_{n-1},\ I(\mu ):=\{k\leq n-1::m_{k}=0\},\ 
l=(l_{1},...,\ l_{n-1})\in {\mathbb{N}}^{n-1}$\ \par
and\ \par
\quad \quad \quad \quad $\displaystyle J^{(0)}_{\mu ,l}(\varphi )(k):=\frac{1}{\pi ^{n}}\displaystyle \int_{{\mathbb{C}}^{n}}^{}{\varphi 
(\zeta )\zeta _{n}^{k}\displaystyle \prod_{i\in I(\mu )}^{}{\zeta _{i}^{l_{i}}\displaystyle 
\prod_{j\notin I(\mu )}^{}\mathds{1}_{C^{(j)}_{m_j}(z,\zeta)}(z_j)\frac{(z _{j}-c_{m_{j},j}(z,\zeta))^{l_j+1}}{(\zeta _{j}-c_{m_{j},j}(z,\zeta))^{-l_{j}-1}}}\,dm_{n}(\zeta 
)}$\ \par
\quad  $\displaystyle J^{(j)}_{\mu ,l}(\varphi )(k):=\frac{1}{\pi ^{n}}\displaystyle \int_{{\mathbb{C}}^{n}}^{}{\varphi 
(\zeta )(\zeta _{n}-c_{j,n})^{-k-1}\displaystyle \prod_{i\in I(\mu )}^{}{\zeta 
_{i}^{l_{i}}\displaystyle \prod_{s\notin I(\mu )}^{}\mathds{1}_{C^{(s)}_{m_s}(z,\zeta)}(z_s)\frac{(z _{s}-c_{m_{s},s}(z,\zeta))^{l_s+1}}{(\zeta _{s}-c_{m_{s},s}(z,\zeta))^{-l_{s}-1}}}\,dm_{n}(\zeta 
)};$\ \par
where 
\begin{eqnarray*}c_{h,k}(z,\zeta)&=&c_{h,k}(z_1,\ldots, z_{k-1},\zeta_{k+1},\ldots,\zeta_n)\qquad 1<k<n\\
c_{h,1}(z,\zeta)&=&c_{h,1}(\zeta_2,\ldots, \zeta_n)\\
c_{h,n}(z,\zeta)&=&c_{h,n}(z_1,\ldots, z_{n-1})\;.\end{eqnarray*}
and the same notation is used for $\mathds{1}_{C^{(j)}_k(z,\zeta)}(z_j)$. We have the link :\ \par
\begin{Theorem}\label{teo_suff_cond} If $J^{(0)}_{\mu ,l}(\varphi )=0$ for every $\mu \in M_{n-1}$ 
and $l\in {\mathbb{N}}^{n},$ then $[\varphi _{n}]_{n}=0\ ;$ given also 
$j=1,...,N_{n},$ if $J^{(j)}_{\mu ,l}(\varphi )=0$ for every $\mu \in M_{n-1}$ 
and $l\in {\mathbb{N}}^{n},$ then $[\varphi _{n},j]_{n}=0.$\ \par
\end{Theorem}
\noindent{\bf Proof: } By direct calculation, using the series expansions given by Lemma \ref{lmm_series_exp}, we have that
$$[K^{(h)}_0(\psi)]_{h+1}(k)=\frac{1}{\pi}A_0^{(h)}\mathds{1}_{C_0^{(h)}}(z_h)\sum_{l\geq 0}z_h^{-l}\int_{\C}[\psi]_h(l)\zeta_{h+1}^kdm_1(\zeta_{h+1})$$
$$[K^{(h)}_0(\psi),m]_{h+1}(k)=\frac{1}{\pi}A_0^{(h)}\mathds{1}_{C_0^{(h)}}(z_h)\sum_{l\geq 0}z_h^{-l}\int_{\C}\frac{[\psi]_h(l)}{(\zeta_{h+1}-c_{m,h+1})^{k+1}}dm_1(\zeta_{h+1})$$
$$[K^{(h)}_j(\psi)]_{h+1}(k)=\frac{1}{\pi}A_j^{(h)}\mathds{1}_{C_j^{(h)}}(z_h)\sum_{l\geq 0}(z_h-c_{j,h})^{l+1}\int_{\C}[\psi,j]_h(l)\zeta_{h+1}^kdm_1(\zeta_{h+1})$$
$$[K^{(h)}_j(\psi),m]_{h+1}(k)=\frac{1}{\pi}A_j^{(h)}\mathds{1}_{C_j^{(h)}}(z_h)\sum_{l\geq 0}(z_h-c_{j,h})^{l+1}\int_{\C}\frac{[\psi,j]_h(l)}{(\zeta_{h+1}-c_{m,h+1})^{k+1}}dm_1(\zeta_{h+1})\;.$$
Therefore, by induction, we obtain that
$$[K^{(n-1)}_{\mu_{n-1}}\cdots K^{(1)}_{\mu_{1}}\phi]_n(l_n)=$$
$$\prod_{i=1}^{n-1} A_{\mu_{i}}^{(i)}\prod_{i\in I(\mu)}\mathds{1}_{C_{\mu_{i}}^{(i)}}(z_i)\sum_{l'\in\N^{n-1}}\prod_{i\in I(\mu)}z_i^{-l_i} J_{\mu,l'\cup\{l_n\}}^{(0)}(\phi)$$
and
$$[K^{(n-1)}_{\mu_{n-1}}\cdots K^{(1)}_{\mu_{1}}\phi, j]_n(l_n)=$$
$$\prod_{i=1}^{n-1} A_{\mu_{i}}^{(i)}\prod_{i\in I(\mu)}\mathds{1}_{C_{\mu_{i}}^{(i)}}(z_i)\sum_{l'\in\N^{n-1}}\prod_{i\in I(\mu)}z_i^{-l_i} J_{\mu,l'\cup\{l_n\}}^{(j)}(\phi)\;.$$
So, if  $J_{\mu, l}^{(0)}(\phi)=J_{\mu, l}^{(j)}(\phi)=0$, all the coefficients vanish, then
$$[\phi_n]_n(k)=0\qquad[\phi_n,j]_n(k)=0$$
as we wanted. \Box
\ \par
\ \par
\ \par
\begin{Definition} We shall say that $\varphi \in L^{r}_{c}({\mathbb{D}^n}\backslash Z)$ 
verifies the structure conditions if  $J^{(0)}_{\mu ,l}(\varphi )=0$ for 
every $\mu \in M_{n-1}$ and $l\in {\mathbb{N}}^{n},$ and if $J^{(j)}_{\mu ,l}(\varphi )=0$ 
for every $\mu \in M_{n-1}$ and $l\in {\mathbb{N}}^{n}.$\ \par
\end{Definition}
\ \par

\section{The polydisc - $q=n$}

As for now, we don't have a way to deal with the integrals $J^{(m)}_{\mu,l}(k)$ on the domain $\P\setminus Z$, so we turn to the much easier case of the polydisc itself. We look first at the problem for $(0,n)$-forms.

Let $\omega$ be a $(0,n)-$form with $L^{r}_c(\P)$ coefficients; we can find a function $\varphi\in L^{r}_c(\P)$ such that
$$\omega=\varphi d\bar{z}_1\wedge\cdots \wedge d\bar{z}_n\;.$$

In this case, the operators $K^{(m)}$ coincide with the outer corona components $K^{(m)}_0$, so the obstructions to a solution of compact support are given by the integrals $J^{(0)}_{0, l}(k)$, where the subscript $0$ stands for a multi-index of the appropriate length containing only $0$s. We have the following result.
\begin{Lemma} \label{lmm_structure_poly}If there is a current $T$ with compact support in ${\mathbb{D}^n}$ 
such that $\overline{\partial }T=\omega ,$ then we have\ \par
\quad \quad \quad \quad $\forall l\in {\mathbb{N}}^{n-1},\ \forall j=1,...,N_{n},\ 
J^{(0)}_{0 ,l}(\varphi )=0,$\ \par
i.e. $\varphi $ verifies the structure conditions for the polydisc.\ \par
\end{Lemma}
\noindent{\bf Proof: } Let $\{\rho_\epsilon\}\subset\Di(\C^n)$ be a family of functions such that $\rho_\epsilon\to\delta_0$, when $\epsilon\to0$, in the sense of distributions, with $\supp\rho_\epsilon\subset\{|z|<\epsilon\}$ and $\|\rho_\epsilon\|_1=1$.

We write
$$T=T_1d\hat{\bar{z}}_1+\ldots+T_nd\hat{\bar{z}}_n$$
so we have
$$\varphi=\debar_1 T_1+\ldots+\debar_n T_n=t_1+\ldots+t_n$$
where, obviously, every $t_h$ is compactly supported in $\P$.

\medskip

We set $T^\epsilon_h=T_h\star\rho_\epsilon\in \Di(\C^n)$; by standard theorems on convolution, 
$$\supp(T^\epsilon_h)\subseteq\{z\ \vert\ \mathrm{dist}(z,\supp T_h)\leq\epsilon\}$$
so, for $\epsilon$ small enough, all the regularized functions are compactly supported in $\P$ and
$$\debar_h T^{\epsilon}_h=t_h\star\rho_\epsilon=t^\epsilon_h\;.$$

By Lemma \ref{lmm_1var}, we have that
$$[t^\epsilon_h]_h(k)=0$$
for every $k\in\N$ and $h=1,\ldots, n$.

Moreover, we have that
$$\varphi^\epsilon=\varphi\star\rho_\epsilon=t_1^\epsilon+\ldots+t_n^\epsilon$$
and $\varphi^\epsilon\to \varphi$ in $L^{r}$ as $\epsilon\to0$.

\medskip

As $\varphi$ and $\varphi^\epsilon$ are compactly supported in $\P$, for $\epsilon$ small enough, we can see them as continuous functionals on $L^q_{\loc}(\P)$ (where $q^{-1}+r^{-1}
=1$). The convergence $\varphi_\epsilon\to\varphi$ holds also in this sense.

The functions
$\zeta_n^k\prod_{i=1}^n\zeta_i^{l_i}$
are in $L^q_{\loc}(\P)$ for every $l\in\N^{n-1}$, $k\in\N$; therefore
$$J_{0,l}^{(0)}(\phi^\epsilon)(k)\xrightarrow[\epsilon\to0]{}J_{0,l}^{(0)}(\phi)(k)\;.$$

Now, consider $t^\epsilon_h$, with $h\leq n-1$; we know that $[t^\epsilon_h]_h(l)=0$, for every $l$ so we can apply Fubini and get
$$J^{(0)}_{0,l}(t^\epsilon_h)(k)=\frac{1}{\pi^n}\int_{\C^n}t^\epsilon_h(\zeta)\zeta_n^k\prod_{i=1}^n\zeta_i^{l_i}dm_n(\zeta)=$$
$$\frac{1}{\pi^n}\int_{\C^{n-1}}\zeta_n^k\prod_{\substack{i=1\\i\neq h}}^n\zeta_i^{l_i}\int_\C t^\epsilon_h(\zeta)\zeta_h^{l_h}dm_1(\zeta_h)dm_{n-1}(\hat{\zeta}_h)=$$
$$\frac{1}{\pi^n}\int_{\C^{n-1}}\zeta_n^k\prod_{\substack{i=1\\i\neq h}}^n\zeta_i^{l_i}[t^\epsilon_h]_h(l_h)dm_{n-1}(\hat{\zeta}_h)=0\;;$$

If $h=n$, it is again an application of Fubini's theorem to show that $J^{(0)}_{0, l}(t^\epsilon_n)(k)=0$.

By additivity of the integral, it follows that $J^{(0)}_{0, l}(\varphi^\epsilon)(k)=0$, so letting $\epsilon\to0$ we obtain the thesis. \Box

\begin{Theorem} \label{teo_sol_n_poly} If $\omega $ is a $(0,n)$-form 
in $L^{r}_{c}({\mathbb{D}^n})$ such that there is a $(0,n-1)$ 
current $T,$ compactly supported in ${\mathbb{D}^n},$ such that 
$\overline{\partial }T=\omega ,$ then we can find a $(p,n-1)$-form $\eta \in L^{r}_{c}({\mathbb{D}^n})$ 
such that $\overline{\partial }\eta =\omega .$\ \par
\end{Theorem}
\textbf{Proof: }By Corollary \ref{cor_structure} we can write $\varphi=\varphi_1+\ldots+\varphi_n$ and, by Lemma \ref{lmm_1var}, the convolutions
$$f_1=\varphi_1\star_1\frac{1}{\pi z_1},\ \ldots,\ f_{n-1}=\varphi_{n-1}\star_{n-1}\frac{1}{\pi z_{n-1}}$$
are compactly supported and
$$\debar_1 f_1+\ldots+\debar_{n-1} f_{n-1}=\varphi_1+\ldots+\varphi_{n-1}=\varphi-\varphi_{n}\;.$$
Moreover, by Lemma \ref{lmm_structure_poly}, $\phi$ satisfies the structure conditions, then, by Theorem \ref{teo_suff_cond}, $[\phi_n]_n(k)=0$ for every $k\in\N$. So, also 
$$f_n=\varphi_n\star_n\frac{1}{\pi z_n}$$
is compactly supported, always by Lemma \ref{lmm_1var}.

We set
$$\eta=\sum_{j=1}^n(-1)^{j-1}f_jd\hat{\bar{z}}_j$$
so that
$$\debar \eta=\varphi d\bar{z}$$
and the coefficients of $\eta$ belong to $L^{r}_c(\P)$. \Box
\ \par
\begin{Remark}
We have that $\|f_j\|_r\leq\gamma\|\phi_j\|_r$, where $\gamma$ depens only on the dimension $n$ and on the radii of $\P$. We recall that $\|K^{(m)}_0\phi\|_r\leq A_0M\|\phi\|_r$, so $\|f_j\|_r\leq (A_0M+1)^j\gamma\|\phi\|_r$; this means that the linear operator associating to $\omega$ the solution $\eta$ is linear and bounded from $L^{r}_c$ to $L^{r}_c$.
\end{Remark}
\ \par

\section{The polydisc - $q=n-1$}

We turn our attention to $(0,n-1)-$forms. Firstly, we give a refined version of Lemma \ref{lmm_structure_poly}.

\begin{Proposition}\label{prp_struct_adv_poly}Suppose $\varphi\in L^{r}(\C^n)$ and $T_1, \ldots, T_{n-1}$ are distributions, compactly supported in $\P$, such that
$$\varphi=\debar_1T_1+\ldots+\debar_{n-1}T_{n-1}\;.$$
Then we can find $\varphi_1,\ldots, \varphi_{n-1}\in L^{r}(\C^n)$, compactly supported in $P $ such that $\varphi=\varphi_1+\ldots+\varphi_{n-1}$ and $[\varphi_i]_i(k)0$, for every $k\in\N$.\end{Proposition}
\noindent{\bf Proof: } After performing the same regularization as in the proof of Lemma \ref{lmm_structure_poly}, we have
$$\frac{1}{\pi^n}\int_{\C^n}t^\epsilon_h(\zeta)a(\zeta_n)\prod_{i=1}^{n-1}\zeta_i^{l_i}dm_n(\zeta)=0$$
for every $a(\zeta_n)$ for which the integral is well-defined (e.g. $a\in L^1$). This is because $h$ ranges from $1$ to $n-1$, so we can isolate the terms $[t^\epsilon_h]_h(l)$ employing only the functions which appear in the product.

Therefore, the function
$$\frac{1}{\pi^n}\int_{\C^{n-1}}t^\epsilon_h(\zeta)\prod_{i=1}^{n-1}\zeta_i^{l_i}dm_n(\hat{\zeta}_n)$$
vanishes for a.e. $z_n$ and the same is true for the function
$$\frac{1}{\pi^n}\int_{\C^{n-1}}\varphi^\epsilon(\zeta)\prod_{i=1}^{n-1}\zeta_i^{l_i}dm_n(\hat{\zeta}_n)$$
and, letting $\epsilon\to0$, also for
$$\frac{1}{\pi^n}\int_{\C^{n-1}}\varphi(\zeta)\prod_{i=1}^{n-1}\zeta_i^{l_i}dm_n(\hat{\zeta}_n)\;.$$

By the analogue of Theorem \ref{teo_suff_cond} in the first $n-1$ coordinates, 
$$[K^{(n-2)}\cdots K^{(1)}\varphi]_{n-1}(k)=0\;,$$ 
so defining $\varphi_1,\ldots, \varphi_{n-2}$ as in Corollary \ref{cor_structure} and setting $\varphi_{n-1}=\varphi-\varphi_1-\ldots-\varphi_{n-2}$ we have that $[\varphi_i]_i(k)=0$, as requested. \Box

\medskip

The following corollary is immediate.

\begin{Corollary}\label{cor_sol_adv_poly}Given $\omega$ as before and a current $T$, compactly supported in $\P$ such that $\debar T=\omega$, with 
$$T=T_1d\hat{\bar{z}}_1+\ldots+T_{n-1}d\hat{\bar{z}}_{n-1}\;,$$
we can find $\eta$ with $L^{r}(\C^n)$ coefficients, compactly supported in $\P$, such that $\debar \eta=\omega$ and with
$$\eta=\eta_1d\hat{\bar{z}}_1+\ldots+\eta_{n-1}d\hat{\bar{z}}_{n-1}\;.$$
\end{Corollary}

\begin{Remark} Obviously, we can suppose that the coefficient of $d\hat{\bar{z}}_k$ in $T$ is zero and obtain that there exists a solution with coefficients in $L^{r}(\C^n)$ with compact support in $\P $ where the coefficient of $d\hat{\bar{z}}_k$ is zero. 

By induction, we can show that if there exists a solution with the coefficients of $d\hat{\bar{z}}_{k_1},\ldots, d\hat{\bar{z}}_{k_r}$ equal to zero, then we can produce a solution in $L^{r}$ with the same vanishing coefficients.\end{Remark}

We note that the construction of $\varphi_1,\ldots, \varphi_{n-1}$ doesn't involve the $n-$th coordinate, so,  $\debar_n\varphi$ and $\debar_n\varphi_j$ share the same regularity, whatever it is.

\begin{Theorem}\label{teo_sol_n1_poly}If $\omega $ is a $(0,n-1)$-form 
in $L^{r}_{c}({\mathbb{D}^n}),\ \overline{\partial }\omega =0,$ such that $\debar_n\omega_n\in L^{r}, $ then we can find a $(0,n-2)$-form $\beta \in L^{r}_{c}({\mathbb{D}^n})$ 
such that $\overline{\partial }\beta =\omega .$\ \par
\end{Theorem}
{\bf Proof: }We proceed by induction on $n$; the case $n=2$ is true.
If there exists a distribution $T$ with compact support such that $\debar T=\omega$, then, by Corollary \ref{cor_sol_adv_poly}, we have
$$\omega_n=\sum_{j=1}^{n-1}\omega_{nj}$$
with $\omega_{nj}\in L^{r}$ and $[\omega_{nj}]_j(k)=0$.

We consider the following family of compactly supported $(0,n-2)-$forms in $\C^{n-1}$, depending on the parameter $z_n$:
$$\psi_{z_n}=\sum_{j=1}^{n-1}\left(\omega_j+(-1)^{n+j}\frac{\de \omega_{nj}}{\de\bar{z}_n}\star_j\frac{1}{\pi z_j}\right)d\hat{\bar{z_j}}\;.$$
Note that, as $\psi_{z_n}$ is thought as a form in $\C^{n-1}$, the notation $d\hat{\bar{z}}_j$ has to be understood as the exterior product of the differentials $d\bar{z}_1,\ldots, d\bar{z}_{n-1}$, with $d\bar{z}_j$ missing.

Now, we have that
$$(\debar' \psi_{z_n})\wedge d\bar{z}_n=\debar\omega=0$$
where $\debar'$ operates in the first $n-1$ coordinates. We note that
$$\frac{\de}{\de \bar{z}_j}\left(\omega_j+(-1)^{n+j}\frac{\de \omega_{nj}}{\de\bar{z}_n}\star_j\frac{1}{\pi z_j}\right)=\debar_j\omega_j+(-1)^{n+j}\debar_n\omega_{nj}$$
belongs to $L^{r}(\C^n)$ for almost all $z_n$. By inductive hypothesis, we can solve $\debar' \xi_{z_n}=\psi_{z_n}$ with compact support (and the result will be in $L^{r}(\C^n)$).

We have $\debar (\xi_{z_n}\wedge d\bar{z}_n)=\psi_{z_n}\wedge d\bar{z}_n$; we define a $(0,n-2)-$form in $\C^{n}$ with
$$\gamma=\sum_{j=1}^{n-1}(-1)^{j-1}\omega_{nj}\star_{j}\frac{1}{\pi z_j}d\hat{\bar{z}}_{jn}\;.$$
So we have
$$\debar \gamma=\omega_nd\hat{\bar{z}}_n+\sum_{j=1}^{n-1}(-1)^{n+j-2}\frac{\de \omega_{nj}}{\de\bar{z}_n}\star_j\frac{1}{\pi z_j}\;;$$
therefore
$$\debar(\gamma+\xi_{z_n}\wedge d\bar{z}_n)=\omega\;.$$
The form $\gamma+\xi_{z_n}\wedge d\bar{z}_n$ has compact support and belongs to $L^{r}(\C^n)$. \Box

\section{The polydisc - $1<q<n-1$}

\setcounter{equation}{0}\ \par

Let $\omega$ be a generic $(0,q)-$form and let us write
$$\omega=\sum_{|J|=n-q}\omega_Jd\hat{\bar{z}}_J\;.$$
We restate here the condition $(\ast)$ given in the introduction
\begin{equation}\label{eq_der_cond}
(\ast)\qquad\debar_{j_{n-q}}\!\cdots\debar_{j_k}\omega_J\in L^{r}(\C^n)\qquad k=1,\ldots, n-q\;,\quad \forall\ |J|=n-q\;.\end{equation}

\begin{Theorem}\label{teo_sol_gen_poly}If $\omega $ is a $(0,q)$-form 
in $L^{r}_{c}(\mathbb{D}^n),\ \overline{\partial }\omega =0,$ fullfilling condition (\ref{eq_der_cond}), 
then we can find a $(0,q-1)$-form $\beta \in L^{r}_{c}({\mathbb{D}^n})$ 
such that $\overline{\partial }\beta =\omega .$\ \par
\end{Theorem}
{\bf Proof: }Following Hörmander \cite[Chapter 2]{Hormander66}, we can write
$$\omega=g\wedge d\bar{z}_n+h$$
where $g$, $h$ do not contain $d\bar{z}_n$.

We can look at $h$ as a family of $(0,q)-$forms in $\C^{n-1}$, depending on the complex parameter $z_n$; similarly, $g$ can be understood as a family of $(0,q-1)-$forms.

We denote by $\debar_{\C^{n-1}}$ the $\debar$ operator in the first $n-1$ variables, that is
$$\debar_{\C^{n-1}}\psi=\sum_{n\not\in I}\sum_{k\not\in I\cup\{n\}}\!\!\! \debar_k \psi_Id\bar{z}_k\wedge d\bar{z}_I\;.$$
If $\psi$ doesn't contain $d\bar{z}_n$, then $\debar' \psi=\debar_{\C^{n-1}}\psi$.

We proceed by induction on the dimension and we prove the following:
\begin{itemize}
\item[$\mathbf{I_n.1}$] the statement of the theorem holds in $\C^n$ and $\beta$ depends linearly on $\omega$;
\item[$\mathbf{I_n.2}$] if the coefficients of $\omega$ depend on a parameter $z_{n+1}\in\C$ in such a way that $\omega, \debar \omega\in L^{r}_c(\C^{n+1})$, then also $\beta, \debar \beta\in L^{r}_c(\C^{n+1})$, where the $\debar$ is intended in $n+1$ variables.
\end{itemize}

We note that $I_2.1$ and $I_2.2$ hold. We assume $I_{n-1}.1$ and $I_{n-1}.2$ to hold.
\ \par
\ \par
\emph{Reduction. } We note that $\debar_{\C^{n-1}}h=0$; therefore, $h$ is a family of $\debar-$closed $(0,q)-$forms in $\C^{n-1}$ depending on the parameter $z_n$. Moreover, by assumption, $\debar_nh_I\in L^{r}_c(\C^{n})$. We denote by $U_{t}$ the $(n-1)-$dimensional open set $\P\cap\{z_n=t\}$ and we note that $U_t$ is still a polydisc, hence Stein, for every $t$ for which it is non-empty.

As a well known consequence of Serre's duality (see \cite{Serre55}) we have $H^{q}_c(U_t,\Ol)=0$, if $2\leq q\leq n-2$; therefore we can find a family $T$ of $(0,q-1)-$currents in $\C^{n-1}$ such that $\debar_{\C^{n-1}} T=h$ for almost every $z_n$. Then, by $I_{n-1}.2$, we can find a family $H$ with $H\in L^{r}_c(\P)$ (and therefore $H_{z_n}\in L^{r}_c(U_{z_n})$ for almost every $z_n$) and with $\debar H\in L^{r}_c(\C^n)$.

Moreover, as $H_{z_n}$ depends linearly on $h$ by $I_{n-1}.1$, if $h_{z_n}=0$, then also $H_{z_n}=0$. Therefore, $H$ is compactly supported in $P$.

Now,
$$\debar H=\debar_{\C^{n-1}}H+ \sum_{I}\debar_n H_Id\bar{z}_n\wedge d\bar{z}_I=h+\sum_{I}\debar_n H_Id\bar{z}_n\wedge d\bar{z}_I$$
so
$$\omega-\debar H=g'\wedge d\bar{z}_n$$
where $g'$ does not contain $d\bar{z}_n$. Moreover, as $\omega$ and $\debar H$ are in $L^{r}_c(\P)$, also $g'$ is. Further, we observe that
$$(\debar_{\C^{n-1}}g')\wedge d\bar{z}_n=\debar(\omega-\debar H)=\debar \omega=0\;$$
and finally, for $z_n$ fixed, $g'$ is a$(0,q-2)-$form in $\C^{n-1}$, fullfilling condition (\ref{eq_der_cond}).

\ \par

\emph{Solution. } We have reduced ourselves to solve $\debar G= g'\wedge d\bar{z}_n$, but as $\debar_{\C^{n-1}}g'=0$, we can, by the same argument used in the reduction, obtain a family $G'$ of $(0,q-2)$ forms in $\C^{n-1}$ such that $\debar_{\C^{n-1}}G'=g'$, by $I_{n-1}.2$.

Again, by the same reasoning, $G'\in L^{r}_c(\P)$ and if we set $G=G'\wedge d\bar{z}_n$, we obtain a $(0,q-1)-$form $G\in L^{r}_c(\P)$ such that $\debar G=g'\wedge d\bar{z}_n$.

So, $\beta=G+H$ is the solution we looked for. This shows $I_n.1$.

To show $I_n.2$ it is enough to notice that all our operations are constructive and preserve the regularity (or summability) of an extra parameter. \Box

\ \par
\begin{Remark} We have to separate the case of $(0,n-1)$-forms from the 
general case because in that case Serre's duality tells us only that $H^{n-1}_c(U_t,\Ol)$ is equal to the topological dual of $H^{0}(U_t,\Omega^{n-1})$, in general not vanishing, so the induction doesn't work there. \ \par
\end{Remark}

\begin{Remark} We note that, in the proof of Theorem \ref{teo_sol_gen_poly}, we never actually used the fact that our domain is the polydisc. Indeed, if we had the analogues of Theorems \ref{teo_sol_n_poly} and \ref{teo_sol_n1_poly} for the domain $\P\setminus Z$ in every dimension, then we could apply the same proof to get Theorem \ref{teo_sol_gen_poly} for $\P\setminus Z$, with exactly the same statement.\end{Remark}

As a corollary of the previous results, we obtain the following.

\begin{Corollary}\label{cor_vanishing}
Let $\omega$ be a $(0,q)-$form with compact support in $\P\setminus Z$ and satisfying conditions \ref{eq_der_cond}, then, for any $k\in\N$, we can find a $(0,q-1)-$form $\beta\in L^r_c(\P)$ such that $\debar(f^k\beta)=\omega$. Equivalently, we can find a $(0,q-1)-$form $\eta=f^k\beta$ such that $\eta\in L^r_c(\P)$, $\eta$ is $0$ on $Z$ up to order $k$ and $\debar \eta=\omega$.\end{Corollary}

\noindent{\bf Proof: } The $(0,q)-$form $\phi:=\omega/f^k$ is still $\debar-$closed and satisfies \ref{eq_der_cond}; hence we have a $(0,q-1)-$form $\beta\in L^r_c(\P)$ such that $\debar \beta=\phi$. So $\eta=f^k\beta$ verifies all the requirements. $\hfill$ \Box

\bibliographystyle{/usr/local/texlive/2009/texmf-dist/bibtex/bst/base/plain}

\begin{thebibliography}{1}

 \bibitem{AbdelkaderKhidr04}
{\sc O.~Abdelkader and S.~Khidr}, {\em Solutions to
  {$\overline\partial$}-equations on strongly pseudo-convex domains with
  {$L^p$}-estimates}, Electron. J. Differential Equations,  (2004) (electronic).
  
\bibitem{AmarMatheron04}
{\sc E.~Amar and E.~Matheron},
\newblock {\em Analyse Complexe},
\newblock Cassini, 2004.

\bibitem{AndreottiGrauert62}
{\sc A.~Andreotti and G.~Grauert}, {\em Th\'eor\`emes de finitude pour la cohomologie des espaces complexes}, Bull. Soc. Math. France, 90, 1962, pp.~193--259.
 
  \bibitem{ChangLee00}
{\sc C.~H. Chang and H.~P. Lee}, {\em {$L^p$} estimates for
  {$\overline\partial$} in some weakly pseudoconvex domains in {${\bf C}^n$}},
  Math. Z., 235 (2000), pp.~379--404.
  
 \bibitem{Fischer01}
{\sc B.~Fischer}, {\em
$L^p$ estimates on convex domains of finite type}, 
Math. Z. 236 (2001), no. 2, pp.~401--418.  
  
 \bibitem{fornsib}
{\sc J.~E. Fornaess and N.~Sibony}, {\em On $\it{L}^p$ estimates for
  $\overline{\partial}$}, in In Several Complex Variables and Complex Geometry,
  Part 3, American Mathematical Society, 1989, pp.~129--163.
  
  \bibitem{Hormander66}
  {\sc L. H\"ormander,} {\em
An introduction to complex analysis in several variables},  D. Van Nostrand Co., Inc., Princeton, N.J.-Toronto, Ont.-London 1966
  
 \bibitem{Jouenne00}
   {\sc C.~Jouenne}
    {\em Convexit\'e et \'equations de {C}auchy-{R}iemann avec
              estimations {$L^p$}}, Publ. Mat., {44}, {2000}, {1}, pp.~{309--323}
  
 \bibitem{Kerzman70}
{\sc N.~Kerzman}, {\em H\"older and {$L^{r}$} estimates for solutions of {$\bar
  \partial u=f$} in strongly pseudoconvex domains}, Bull. Amer. Math. Soc., 76
  (1970), pp.~860--864.

\bibitem{Kerzman71}
\leavevmode\vrule height 2pt depth -1.6pt width 23pt, {\em H\"older and
  {$L^{r}$} estimates for solutions of {$\bar \partial u=f$} in strongly
  pseudoconvex domains}, Comm. Pure Appl. Math., 24 (1971), pp.~301--379.

\bibitem{Khidr08}
{\sc S.~Khidr}, {\em Solving {$\overline{\partial}$} with {$L^p$}-estimates on
  {$q$}-convex intersections in complex manifold}, Complex Var. Elliptic Equ.,
  53 (2008), pp.~253--263. 

\bibitem{Krantz76}
{\sc S.G. ~Krantz}, {\em
Optimal Lipschitz and $L^{p}$ regularity for the equation $\overline \partial u=f$ on stongly pseudo-convex domains},  
Math. Ann. 219 (1976), no. 3, pp.~233--260. 

\bibitem{Landucci79}
{\sc M.~Landucci}, {\em Solutions with ``precise'' compact support of the
  {$\bar \partial $}-problem in strictly pseudoconvex domains and some
  consequences}, Atti Accad. Naz. Lincei Rend. Cl. Sci. Fis. Mat. Natur. (8),
  67 (1979), pp.~81--86 (1980).

\bibitem{Landucci80}
\leavevmode\vrule height 2pt depth -1.6pt width 23pt, {\em Solutions with
  precise compact support of {$\bar \partial u=f$}}, Bull. Sci. Math. (2), 104
  (1980), pp.~273--299.
  
  \bibitem{Laurent11}
  {\sc C.~Laurent-Thi\'ebaut},  {\em Holomorphic Function Theory in Several Variables}, Springer London, 2011.
  
 \bibitem{Li10}
{\sc X.-D. Li}, {\em {$L^p$}-estimates and existence theorems for the
  {$\overline\partial$}-operator on complete {K}\"ahler manifolds}, Adv. Math.,
  224 (2010), pp.~620--647.
  
\bibitem{Serre55}
{\sc J.-P. Serre}, {\em Un th\'eor\`eme de dualit\'e}, Comment. Math. Helv., 29
  (1955), pp.~9--26.
  
   \bibitem{Shapiro09}
  {\sc A.~Shapiro, D.~Dentcheva, A.P.~Ruszczy\'nski},  {\em Lectures on stochastic programming: modeling and theory}, SIAM, 2009
  
\end{thebibliography}

\end{document}